\documentclass[12pt]{article}
\title{Classification of $osp(2|2)$ Lie super-bialgebras}
\author{Cezary Juszczak \\
\normalsize Institute for Theoretical Physics,
\normalsize University of Wroc\l aw\\ 
\normalsize pl.\ M.\ Borna 9, Wroc\l aw,
\normalsize Poland}

\def\k#1{\ref{#1}}

\def\ztwelve{z_{12}}
\def\zone{z_{1}}
\def\ztwo{z_2}
\def\zthree{z_3}
\def\zfour{z_4}
\def\zfive{z_5}
\def\zsix{z_6}
\def\zseven{z_7}
\def\zeight{z_8}
\def\znine{z_9}
\def\zten{z_{10}}

\def\zonethree{z_{13}}
\def\zonefour{z_{14}}
\def\zonefive{z_{15}}
\def\zonesix{z_{16}}
\def\ztwoone{z_{21}}
\def\ztwotwo{z_{22}}
\def\zzten{Z_{10}}
\def\zzonesix{Z_{16}}
\def\zzonefive{Z_{15}}

\def\fr{\frac}

\newtheorem{lemma}{Statement}
\newtheorem{Def}{Definition}
\def\blem{\begin{lemma}}
\def\elem{\end{lemma}}
\def\bel{\begin{equation}\label}
\def\ee{\end{equation}}
\def\ba{\begin{array}}
\def\ea{\end{array}}
\def\ars{\\[2mm]  }
\def\dsp{\displaystyle}
\def\r#1{(\ref{#1})}
\def\tens{\otimes}
\def\t{\tilde }
\def\rank{\mathop{\rm rank}}
\def\diag{\mathop{\rm diag}}
\def\ben{\begin{enumerate}}
\def\een{\end{enumerate}}

\def\beq{\begin{eqnarray}}
\def\eeq{\end{eqnarray}}
\def\c#1#2#3{c_{#1 #2}{}^{#3}}
\def\f#1#2#3{f_{#1}{}^{#2#3}} 
\def\d{\delta}
\date{}

\begin{document}
\maketitle
\begin{abstract}
The co-Lie structures compatible with the $osp(2|2)$ Lie super algebra structure are
investigated and found to be all of coboundary type. The corresponding
classical $r$-matrices are classified into several disjoint families.
The $osp(1|2) \oplus u(1)$ Lie super-bialgebras are also classified.
\end{abstract}

\section*{Introduction}
The need of classification of Lie bialgebras \cite{drinfeld} comes from their close
relation with $q$-deformations of universal enveloping 
algebras in the Drinfeld sense.
 To each such deformation there 
corresponds a Lie bialgebra which may be recovered from the first order of the 
deformation of the coproduct.

It has also been shown \cite{eting} that each Lie bialgebra admits 
quantization. So the classification of Lie bialgebras 
can be seen as the first step in classification of quantum algebras.

Along these lines several efforts
(see e.g.\ \cite{zakrz,b,c,balles} to list only a few)
have been undertaken in order to classify those Hopf algebras which 
can be of importance in physics. 

The $Osp(2|2)$ super-group is a subgroup of two-dimensional $N=2$ 
super-conformal symmetry which plays an important role in string theory. 
In \cite{conform} the correlation functions of $N=2$ super-conformal 
field theory were found by using the $Osp(2|2)$ symmetry group.
Lattice models based on $U_q(osp(2|2))$ symmetry were constructed in 
\cite{lattice}, where also new solutions to the graded Yang-Baxter
equation were found.
 
A few examples of quantum deformations of $osp(2|2)$
\cite{d, preeti, deguchi, aizawa} were given so far and it 
became evident that their classification would be of much value. 

In this paper we perform a complete classification 
of Lie super-bialgebras $osp(2|2)$ based
on the brut-force computer approach combined with careful
identification of equivalent structures. 
We also classify of the $u(1) \oplus osp(1|2)$ 
Lie super-bialgebras. $u(1) \oplus osp(1|2)$ is the simplest 
central extension of the $osp(1|2)$ subalgebra and is similar 
to $osp(2|2)$ in the fact that it containes $gl(2)$ and $osp(1|2)$ as
subalgebras.
In both cases all the obtained structures are coboundary, 
allowing for a brief exposition of the results in the form 
of list of classical $r$-matrices.

\section{Lie super-bialgebra $osp(2|2)$ }
The $osp(2|2)$ Lie superalgebra $G=G_0 \oplus G_1$
is spanned by the generators ($g_1$, \ldots, $g_8$) $=$
($H$, $X_+$, $X_-$, $B$, $V_+$, $V_-$, $W_+$, $W_-$),  
where $H$, $X_\pm$, $B$ span the subspace $G_0$ of grade 0, 
 and $V_\pm$, $W_\pm$ span the subspace $G_1$ of grade 1.
We refer to the elements of $G_0$ and $G_1$ as bosons 
and fermions respectively. The generators fulfill the following relations:
\bel{r1}
\ba{l l}
[H,X_\pm]= \pm X_\pm,\qquad &    [X_+,X_-]= -2H\,, \ars
[H,B]=0, &[X_\pm, B]=0\,, \ars
[H,V_\pm]=\pm \frac12{V_\pm} & [H,W_\pm]=\pm \frac12{W_\pm}\,, \ars
[B,V_\pm]=\frac12{V_\pm} & [B,W_\pm]=-\frac12{W_\pm}\,, \ars
[X_\pm, V_\pm]=0\,, & [X_\pm,W_\pm]=0\,,\ars
[X_\pm, V_\mp]= \mp V_\pm \,, & [X_\pm, W_\mp]= \mp W_\pm\,,\ars
\{V_\pm,V_\pm\} = \{V_\pm,V_\mp\} =&  \{W_\pm,W_\pm\} = \{W_\pm,W_\mp\} = 0\,, \ars
\{V_+,W_-\} = H - B, & \{W_+,V_-\} = H + B\,,\ars
\{V_\pm,W_\pm\} = X_\pm\,.\ars
\ea
\ee

From the last three relations it is evident 
 that the superalgebra is generated by its fermionic sector $G_1$. 

For the approach taken in the  present paper it is most convenient to use 
the definition of Lie super-bialgebra in terms of the structure constants.
Thus we define 
\begin{Def}
 Lie super-bialgebra \cite{andru} is a vector space $G$ with two linear 
mappings:  
\bel{1.1}
[~.~]: ~ G \tens G \ni g_i \tens g_j \mapsto [g_i,g_j] = \c ijk g_k \in G\,,
\ee
\bel{1.2}
\d : ~ G \ni g_i \mapsto \d(g_i) = \f ijk g_j \tens g_k \in G \tens G\,.
\ee
$[~,~]$ is a Lie bracket on $G$ which means that its structure constants 
$\c ijk$  satisfy the relations:
\bel{1.3}
\c ijk= 0 \mbox{ if  $grade(i)+ grade(j) \not \equiv grade(k)~ (mod~2)$  }
\ee
\bel{1.4}
\c ijk = - z(i,j) \c jik
\ee
\bel{1.5}
\c ijk \c klm z(i,l) + \c jlk \c kim z(j,i) + \c lik \c kjm z(l,j) =0
\ee
where
\bel{1.5a} 
 z(i,j) \equiv (-1)^{grade(g_i) \cdot grade(g_j)}
\ee
$\d^*$ defines Lie bracket on the dual space $G^*$ so its 
structure constants $\f ijk$ fulfill similar relations:
\bel{1.6}
\f kij= 0 \mbox{ if  $grade(i)+ grade(j) \not \equiv grade(k)~ (mod~2)$  }
\ee
\bel{1.7}
\f kij = - z(i,j) \f kji
\ee
\bel{1.8}
\f ikj \f jlm z(k,m) + \f ilj \f jmk z(l,k) + \f imj \f jkl z(m,l) =0.
\ee
Moreover, the two mappings need to be compatible:
\bel{1.9}
\c ijk \f klm = \f ilk \c jkm + \c kjl \f ikm z(m,j) +
                \c jkl \f jkm + \f jlk \c ikm z(i,l)
\ee
\end{Def}

If there exists an element $r= r^{ij} g_i \tens g_j \in G\tens G$ such that :
\bel{1.10}
\d(g_i) = [ r , g_i \tens 1 + 1 \tens g_i ] 
\ee
or, in terms of the structure constants,
\bel{1.11}   
\f ijk = r^{jm} \c mik - \c imj r^{mk}\,
\ee
then the  $G$ is called {\em coboundary} Lie super-bialgebra.
It is easy to see that the graded antisymmetric part of $r$ defined by
$\hat r_{ij}=(r_{ij}-z(i,j)r_{ji})/2$ yields the same $\f ijk$ so we will 
assume that $r \in G\wedge G$. 
\bel{1.12}   
r^{ij}= - r^{ji} z(i,j)\,. 
\ee
Similarly, it can be shown that
projection of $r$ on the $G_0\wedge G_1$ subspace of $G\wedge G$ cannot 
influence $\f ijk$ without violating the condition \r{1.6}.
After subtracting it from $r$ we obtain 
even $r$-matrix $r\in G_0\wedge G_0 \oplus G_1\wedge G_1$.
i.e.\
\bel{1.11a}   
r^{ij}=0 ~\mbox{ if $grade(g_i) \neq grade(g_j)$.} 
\ee

We start with the Lie bracket $[~,~]$ given by \r{r1} and look for 
all the co-brackets $\d$ which are compatile with it.

The commutation relations \r{r1} fix the structure constants  $c_{ij}{}^k$. 
Then our task is to find all the 
$f_i{}^{jk}$ that fullfill \r{1.6}, \r{1.7}, \r{1.8}, and \r{1.9}. 
To this end we use a computer and a symbolic algebra program REDUCE.
We use \r{1.6}, \r{1.7} just to reduce the number of unknowns,
then we solve the set of linear equations \r{1.9} coming from the cocycle 
condition. At this point we are able to obtain a 16-parameter family 
of solutions and by solving the relations \r{1.11} 
are able to find the corresponding classical $r$-matrix.

This leads to the conclusion that all the solutions are coboundary.
It is well known \cite{drinfeld} that the Lie bialgebras of simple Lie
algebras are all of coboundary type. Since the $osp(2|2)$ is a simple Lie 
superalgebra the fact that all its bialgebras are 
coboundary can probably 
be justified from the cohomological point of view. 

We substitute the results into the quadratic equations \r{1.8}
representing the co-Jacobi identity. Solving them 
yields 22 solutions, each  parametrized by up to 6 complex numbers.
Substituting these solutions into the generic $r$-matrix 
we obtain 22 families of classical $r$-matrices.

We consider two coalgebra structures $\d$ (and their corresponding $r$-matrices) 
equivalent if they differ only by a linear transformation of the
generators which preserves the algebra commutation relations \r{r1}.

In the next Section we prove that these transformations form a group 
 $GL(2)\oplus Z_2$.
In Section \ref{details} we use this 4-parameter symmetry
to obtain families of {\em nonequivalent} $r$-matrices
parametrized by at most 2 complex numbers.

\section{Automorphisms of the algebra}

We consider two coalgebra structures $\d$ equivalent if they differ only 
by a change of basis. As a change of the basis we allow only
such linear transformations of the generators that:
(a) parity is preserved, and
(b) algebra structure constants $\c ijk$ are unaffected.
(i.e.\ automorphisms of the Lie super-algebra).

Since the algebra  is generated by the fermions, every such transformation
is generated by a transformation within the fermionic sector $G_1$
which in turn can be identified as a nonsingular $4\times 4$ matrix $A_F$:
\bel{r2}
\pmatrix{\tilde V_+ \cr \tilde V_- \cr \tilde W_+ \cr \tilde W_-}
= A_F \cdot \pmatrix{V_+ \cr V_- \cr W_+ \cr W_-}
\ee
such that $\tilde V_+$, $\tilde V_-$, $\tilde W_+$ $\tilde W_-$ fulfill the
relations \r{r1}.

\blem  The matrix $A_F$ must be either block diagonal or
block anti-diagonal i.e.\ it is of the form
\bel{r7}
A_1=\pmatrix{  A_{VV} & 0 \cr 0 & A_{WW} }\mbox{ or }
A_2=\pmatrix{0 &  A_{VW} \cr  A_{WV}  &  0}\,.
\ee
\elem
{\bf Proof:}

Let us assume the following general expression for $\tilde V_+$
\bel{r5}
\tilde V_+ =  a V_+  + b V_-  + c W_+  +  d W_- \,.
\ee
Then
\bel{r6}
\begin{array}{lll}
0&=& \{ \tilde V_+,\tilde V_+ \}/2  \\
 &=&  \{a V_+  + b V_-  + c W_+  +  d W_- ,\, a V_+  + b V_-  + c W_+  +  d W_-  \}/2\\
 &=& ac X_+ + (ad+bc) H  +  (-ad+bc) B + bd X_- \,.
\end{array}
\ee
The condition $ 0= ac = ad+bc  =-ad+bc =  bd $  has two solutions
($a=b=0$ or $c=d=0$) which shows that $\tilde V_+$ is either a
combination od $V$'s or a combination of $W$'s. Similar reasoning is
valid for $\tilde V_-$, $\tilde W_+$, and $\tilde W_-$.

Now we show that both $\t V_+$ and $\t V_-$ belong to
same sector ($V$ or $W$). Indeed,  assumption to the contrary, i.e.\ that
$\tilde V_+ = a V_+ +b V_- $,  $\tilde V_- = c W_+ +d W_- $
would imply:
\bel{prof2}
\begin{array}{lll}
0&=& \{ \tilde V_+ , \tilde V_- \} \\
&=& \{ a V_+ +b V_- ,\,  c W_+ +d W_- \}\\
&=& ac X_+ +(ad+bc) H   +  (-ad+bc) B +  bd X_-\,.
\end{array}
\ee
Then $ 0 = ac = ad+bc  = -ad+bc =  bd$ with the only
two solutions being ($a=b=0$ or $c=d=0$) would mean
that either $\tilde V_+ =0$ or $\tilde V_- =0$. This would contradict our 
assumption that $A_F$ is nonsingular.

Now when we have proved that $\tilde V_+$, $\tilde V_-$ belong to the
same sector we see that  $\t W_+$ and $\t W_-$ must belong to the 
other one for the $A_F$ to be nonsingular.
The conclusion is that the matrix $A_F$ is block diagonal
or block antidiagonal.\framebox{}
\vskip 3mm
{\bf Comment:} Every block antidiagonal matrix $A_2$ can be written
in the form:
$$
A_2= \pmatrix{
0 & 0& 1 & 0 \cr
0 & 0& 0 & 1 \cr
1 & 0& 0 & 0 \cr
0 & 1& 0 & 0 }
 \cdot A_1
$$ 
where $A_1$ is block diagonal.

\blem The diagonal blocks $A_{VV}$ and $A_{WW}$ of $A_1$
are proportional to each other:
\bel{r8}
A_{VV} =  k \cdot A_{WW}
\ee
where $k=\det A_{VV}$.
\elem
{\bf Proof:}

Assume that
\bel{ala1}
A_{VV} =\pmatrix {a &b \cr c &d } \,,\qquad
A_{WW} =\pmatrix {x &y \cr z &t } \,,
\ee
or equivalently:
\bel{ala2}
\ba{l}
\t V_+ = a V_+ + b V_-\,, \qquad \t W_+ = x W_+ + yW_-\,, \\
\t V_- = c V_+ + d V_-\,,  \qquad \t W_- = z W_+ + tW_-\,,.
\ea
\ee
Then
\bel{as1}
\t X_+ = \{ \t V_+ , \t W_+ \}
= ax X_+ + (bx+ay) H + (bx-ay) B + by X_- \,,
\ee
\bel{as2}
\t X_- = \{ \t V_- , \t W_- \}
= cz X_+ + (dz + ct) H + (dz - ct )B + dt X_-\,.
\ee
Inserting \r{ala2}-\r{as2} into
\bel{ola1}
[\t X_+, \t V_-]= - \t V_+ \,, \qquad [\t X_-, \t V_+]= \t V_- \,,
\ee
we obtain
\bel{ola2}
\ba{lcr}
(ad-bc)(-xV_+ -yV_-) &=& -a V_+ - b V_- \,,\\
(ad-bc)(~ zV_+ ~+ ~tV_-) &=& c V_+ + d V_- \,.
\ea
\ee
from which follows
\bel{ola3}
\pmatrix{ a& b \cr c& d} = (ad-bc)
\pmatrix{ x & y \cr z & t}\,.
\ee
\framebox{}
\vskip 3mm
The remaining relations \r{r1} do not lead to further constraints on
the numbers $a$, $b$, $c$ $d$.
Altogether, we have just shown that
\blem The matrix $A_F$ has the following general form
\bel{r30}
A_F=  \pmatrix{
0 & 0& 1 & 0 \cr
0 & 0& 0 & 1 \cr
1 & 0& 0 & 0 \cr
0 & 1& 0 & 0 }^m
\pmatrix{
a & b& 0 & 0 \cr
c & d& 0 & 0 \cr
0 & 0& a/k & b/k \cr
0 & 0& c/k & d/k }\,,
\qquad m=0,1
\ee
where $a$, $b$, $c$, $d$ are arbitrary complex numbers such that
\bel{r4}
k \equiv \det\pmatrix{ a & b \cr c& d} \neq 0\,.
\ee
\elem
The action of the above symmetry on bosons is defined by the matrix $A_B$
\bel{bos2}
\pmatrix{\t H \cr \t X_+ \cr \t X_- \cr \t B } =
A_B \cdot
\pmatrix{H \cr  X_+  \cr X_- \cr  B }
\ee
where
\bel{bos}
A_B=\pmatrix{ 1 & 0 & 0 & 0 \cr
      0 & 1 & 0 & 0 \cr
      0 & 0 & 1 & 0 \cr
      0 &  0 &  0 & (-1)^m
}
k^{-1} \pmatrix{ ac+ bc & ac & bd & 0 \cr
      2ab & a^2 & b^2 & 0 \cr
      2cd & c^2 & d^2 & 0 \cr
      0 &  0 &  0 & k
}
\ee

We define block diagonal matrix $A={\rm diag}(A_B,A_F)$ such that
\bel{A}
\t {\vec  g} = A \cdot \vec g
\ee
where $\vec g ^T =(g_1,\ldots\, ,g_8)=(H,X_+,X_-,B,V_+,V_-,W_+,W_-)$.
Since all the obtained coalgebra structures are coboundary,
which means that the co-Lie bracket $\delta$ is defined in terms
of classical $r$-matrix:
\bel{rmat1}
\delta(x) = [ r, x \tens 1 + 1 \tens x ] \,,
\ee
it is useful to know the action of the symmetry $A$ on 
classical $r$-matrices.
\bel{rmat2}
r= r^{ij} g_i \tens g_j  \,.
\ee
From
\bel{rmat21}
r= r^{ij} g_i \tens g_j   =
\t r^{ij}\t g_i \tens\t g_j =
\t r^{kl} A_k{}^i g_i   \tens A_l{}^j g_j =
\t r^{kl} A_k{}^i A_l{}^j  g_i   \tens  g_j\,,
\ee
we see that
\bel{sym2}
 r = A^T \t r A  \qquad \mbox{ and } \qquad
\t r = (A^{-1})^T r (A^{-1})\,.
\ee

Because we require the $r$-matrix to be even  \r{rmat2}
is equivalent to
\bel{rmat3}
r= \sum_{i,j=1}^4 r_{ij} g_i \tens g_j +\sum_{i,j=5}^8 r_{ij} g_i \tens g_j =
r_B + r_F\,,
\ee
where we define
\bel{rmat4}
(r_B)_{ij}=r_{ij}\,, \qquad (r_F)_{ij}=r_{i+4\,j+4}\,, \qquad \mbox{ for $i,j=1,2,3,4$\,.}
\ee
Then it follows from \r{sym2}  that:
\bel{sym3}
\t r_F = (A_F^{-1})^T r_F (A_F^{-1})\,,\qquad \t r_B = (A_B^{-1})^T r_B (A_B^{-1}) \,.
\ee

We will use this symmetry to identify which classes of $r$-matrices 
differ only by a change of basis. Furthermore, in most cases we will be able
to eliminate some parameters from $r$ by an appropriate choice of
parameters $a$, $b$, $c$, $d$ of $A_F$.

$$
r_F= \pmatrix{ r_{VV} & r_{VW} \cr r_{WV} & r_{WW} }
$$
where $r_{WV}=r_{VW}^T$, $r_{VV}=r_{VV}^T$ and $r_{WW}=r_{WW}^T$. 

The action of the symmetry $S$  consists in exchanging 
$r_{VV}$ with $r_{WW}$ and $r_{WV}$ with $r_{VW}$ or, in other words, 
exchanging  $V$ with $W$.

The action of symmetry $A$ on $r_{WW}$ and $r_{VV}$ is the following
$\tilde r_{VV}= A^T r_{VV} A $, $\tilde r_{WW}= A^T r_{WW} A  (\det A)^{-2}$.
Since they are symmetric we may use the Sylvester theorem to prove that 
for example $r_{VV}$ can be made equal $\diag(1,1)$, $\diag(1,0)$ or 
0 depending on its rank. Once one of these forms is achieved we may use
the remaing symmetry to simplify $r_{WV}$ and then the remaining to simplify 
$r_{WW}$. 

\section{Details of the equivalence considerations \label{details}}

Our strategy of bringing the solutions to the 
'canonical' form by using some change of basis 
can be summarized in the following steps:
\ben
\item If $r_{VV}$ has lower rank than $r_{WW}$ we apply the symmetry $S$ which makes them 
interchange.  Now $\rank r_{VV} \geq \rank r_{WW}$\,. 
It also turns out that $\rank r_{WW}<2 $ after this step. 
\item We can diagonalize $r_{WW}$ so that we obtain 
$$
r_{WW}=\pmatrix{1&0\cr 0 &0} \mbox{ or }r_{WW}=\pmatrix{   0&0 \cr 0 & 0}
$$
depending on the initial rank of $r_{WW}$.
\item Now we try to simplify $r_{WV}$ while preserving the form of $r_{WW}$.
 Depending on whether $r_{WW}$ vanishes or not we have either 
 the full $GL(2)$ symmetry at our disposal or just $A$ of the form 
\bel{as}
A=\pmatrix{1 & b \cr 0 & d} 
\ee

In general the matrix $r_{VW}$ does not have to be symmetric nor antisymmetric.
However, it almost always is. 
Three general possibilities occur in this case:
\ben
  \item When $r_{VW}$ is antisymmetric it is invariant with respect to $A$. 
        We just proceed to simplify $r_{VV}$.
  \item When $r_{VW}$ is symmetric we make sure it is antidiagonal i.e.\ of the 
       form $\pmatrix{ 0 & z \cr z & 0}$ or 
       diagonal $\pmatrix{ x & 0 \cr 0 & 0}$.
  \item In other cases we easily obtain $r_{VW} = \pmatrix{0 & x \cr 0 & 0}$  
\een          
\item 
 We use the remaining symmetry to simplify the part $r_{VV}$. Here what 
 we can achieve depends on the previous steps and on the rank of 
 $r_{VV}$.   
\item We can use the WV scaling ($A_{VV} =\pmatrix{ x & 0 \cr 0 &x }$) to scale $r_{VV}$ with respect to $r_{WW}$.
\item If we still have any freedom  we try to make simpler the bosonic part.
  The main rule is to get rid of $X_-$ if possible.
  
\een

\def\koniec{\\}
\def\pocz{Case }
\def\nad{\over}

The list of computer generated solutions for the 
classical $r$-matrix consists of 22 entries, 
many of them equivalent. We list them below in the form $r_B$, $r_F$. 
The numbering of the cases has been changed from the original computer 
output in order to group equivalent cases. 

 The $r$-matrices which are equivalent  are discused together.
The parameters $a$, $b$, $c$, $d$, and $m$ of $A$  are always given
for the symmetry \r{r30} \r{bos} which makes the $r$-matrices identical.
Here 'identical' means the same up to renaming of some arbitrary constants.

It is well known that computer produced results are the generic ones.
For example if we solve the equation $xy=1$ with respect to $x$ we obtain 
the generic solution $x=\fr1y$ which does not make sense when 
$y=0$. In this case the equation has no solutions. If however we start with 
equation $xy=z$ we obtain $x=\fr z y$. This makes no sense either 
when $y=0$ but when $z=y=0$ there are in reality continuum of solutions. 
They can be recovered from the generic one by taking the limit
$z = \lambda y \mapsto 0$ where $\lambda=$const. 

It is therefore very important to perform the analysis of what happens when 
some parameters of our set of solutions tend to $0$.  Such singular limits were investigated 
when needed and they are listed as footnotes.
We did our best to find out all the special cases which might have been
missing from the computer generated list of results. 

We use the strategy outlined above to discuss each group and come up with 
a list of nonequaivalent $r$-matrices.

\subsection{Cases \k{1},\k{13},\k{5},\k{7},\k{16},\k{19}}
\def\zonefive{K}
\def\zzonefive{K}
\def\zonefour{Y}
\def\zonethree{U}
\def\zonesix{L}
\def\zzonesix{L}
\def\zten{M}
\def\zzten{M}
\def\zfour{J}
\def\zthree{N}
\def\zfive{T}
\def\znine{B}
\def\zeight{C}
\def\zone{X}
\def\ztwo{Z}
\def\ztwelve{F}
\def\zsix{S}
\def\zzsix{S'}
\def\zseven{P}
\def\ztwoone{G}
\def\ztwotwo{H}
\def\ztwotwo{G}

\pocz{\k{1}}
$$
r_{\k{1}} =\pmatrix{0&\zfour\nad 2&( - 2\cdot \zonefive \cdot \zonesix )\nad \zfour&\zonefive  + \zonesix  \cr
      ( - \zfour)\nad 2&0&( - (\zonefive  + \zonesix ))\nad 2&\zfour\nad 2 \cr
 (2\cdot \zonefive \cdot \zonesix )\nad \zfour&(\zonefive  + \zonesix )\nad 2&0&(2\cdot \zonefive \cdot \zonesix )\nad \zfour \cr
  - (\zonefive  + \zonesix ) & ( - \zfour)\nad 2&( - 2\cdot  \zonefive \cdot \zonesix )\nad \zfour&0}\koniec
, $$ $$
\pmatrix{(\zonethree \cdot \zfour)\nad (2\cdot \zonesix )&\zonethree &0&(\zonefive  - \zonesix )\nad 2\cr \zonethree &(2\cdot \zonethree \cdot \zonesix )\nad \zfour&( - \zonefive  + \zonesix )\nad 2&
0\cr 0&( - \zonefive  + \zonesix )\nad 2&0&0\cr (\zonefive  - \zonesix )\nad 2&0&0&0}\koniec
$$
\pocz {\k{13}} 
$$
r_{\k{13}}=\pmatrix{0&( - \zfour)\nad 2&(2\cdot \zten\cdot \zonefive )\nad \zfour&\zten + \zonefive \cr \zfour\nad 2&0&(\zten + \zonefive )\nad 2&\zfour\nad 2\cr ( - 2\cdot
\zten\cdot \zonefive )\nad \zfour&( - (\zten + \zonefive ))\nad 2&0&(2\cdot \zten\cdot \zonefive )\nad \zfour\cr  - (\zten + \zonefive ) & ( - \zfour)\nad 2&( - 2
\cdot \zten\cdot \zonefive )\nad \zfour&0}\koniec
$$
$$
\pmatrix{0&0&0&( - \zten + \zonefive )\nad 2\cr 0&0&(\zten - \zonefive )\nad 2&0\cr 0&(\zten - \zonefive )\nad 2&(\zfour\cdot \zeight)\nad (2\cdot
\zten)&\zeight\cr ( - \zten + \zonefive )\nad 2&0&\zeight&(2\cdot \zten\cdot \zeight)\nad \zfour}\koniec
$$
This is equivalent to case \k{1} by means of the $S$ symmetry 
and then renaming some parameters.\\
\pocz {\k{5}} 
$$
r_{\k{5}}=\pmatrix{0&0&\zthree\nad 2&\zonesix \cr 0&0&( - \zonesix )\nad 2&0\cr ( - \zthree)\nad 2&\zonesix \nad 2&0&( - \zthree)\nad 2\cr  - \zonesix &0&\zthree
\nad 2&0}\koniec ,
\pmatrix{0&0&0&( - \zonesix )\nad 2\cr 0&\znine&\zonesix \nad 2&0\cr 0&\zonesix \nad 2&0&0\cr ( - \zonesix )\nad 2&0&0&0}\koniec
$$
\pocz {\k{7}} 
$$
r_{\k{7}}=\pmatrix{0&0&( - \zthree)\nad 2&\zten\cr 0&0&\zten\nad 2&0\cr \zthree\nad 2&( - \zten)\nad 2&0&( - \zthree)\nad 2\cr  - \zten&0&\zthree
\nad 2&0}\koniec ,
\pmatrix{0&0&0&( - \zten)\nad 2\cr 0&0&\zten\nad 2&0\cr 0&\zten\nad 2&0&0\cr ( - \zten)\nad 2&0&0&\zfive}\koniec
$$
\pocz {\k{16}}
$$
r_{\k{16}}=\pmatrix{0&( - \zfour)\nad 2&0&\zonefive \cr \zfour\nad 2&0&\zonefive \nad 2&\zfour\nad 2\cr 0&( - \zonefive )\nad 2&0&0\cr  - \zonefive &( - \zfour)\nad 2
&0&0}\koniec ,
\pmatrix{0&0&0&\zonefive \nad 2\cr 0&0&( - \zonefive )\nad 2&0\cr 0&( - \zonefive )\nad 2&\zseven&0\cr \zonefive \nad 2&0&0&0}\koniec
$$
\pocz {\k{19}}
$$
r_{\k{19}}=\pmatrix{0&\zfour\nad 2&0&\zonefive \cr ( - \zfour)\nad 2&0&( - \zonefive )\nad 2&\zfour\nad 2\cr 0&\zonefive \nad 2&0&0\cr  - \zonefive &( - \zfour)\nad 2
&0&0}\koniec ,
\pmatrix{\ztwelve&0&0&\zonefive \nad 2\cr 0&0&( - \zonefive )\nad 2&0\cr 0&( - \zonefive )\nad 2&0&0\cr \zonefive \nad 2&0&0&0}\koniec
$$
All the above cases become identical with the case \k{19} after bringing them
to the 'canonical' form.
We take $r_{\k{19}}$ for further consideration.
We have two distinct cases \\
(a) $\zonefive  \neq 0$, \\
(b) $\zonefive  = 0$. \\
In case (a) we apply the symmetry 
$$\pmatrix{a &b \cr c& d} =
\pmatrix{ 1 \nad t & 0 \cr 
-\zfour \nad 2 t \zonefive   & t  }
$$
with $t=\sqrt{\ztwelve}$ if $\ztwelve \neq 0$ and $t=1$ if $\ztwelve =0$. We obtain
$$
\t r_{a1}=\pmatrix{0&0& 0&\zonefive \cr 0&0& - \zonefive \nad 2&0\cr 0&\zonefive \nad 2&0& 0\cr
- \zonefive &0& 0 &0},
\pmatrix{ \alpha &0&0& - \zonefive \nad 2\cr 0& 0&\zonefive \nad 2&0\cr 0&\zonefive \nad 2&0&0\cr
 - \zonefive  \nad 2&0&0&0},
$$
where $\alpha=0$ if $\ztwelve=0$ and $\alpha=1$ otherwise. 
In standard notation we have
\bel{a1}
r_{a1} = x( -2 H \wedge B + X_+ \wedge X_- + V_+ \wedge W_- - V_- \wedge W_+)
+ \alpha \frac12  V_+ \wedge V_+ 
\ee
where $x=-\zonefive /2$. 
In case (b)  we can use the scalings $(+-)$ and $(WV)$ to obtain a nonstandard 
$r$-matrix
\bel{a2}
r_{a2} = \alpha (H - B) \wedge X_+ + \beta V_+ \wedge V_+\,, 
\ee
where $\alpha$, $\beta$ = $0,1$. 

\subsection{Case \k{3}}
\pocz {\k{3}}
$$
r_{\k{3}}=\pmatrix{0&\zonefour &(\zonesix \cdot ( - 2\cdot \zonefive  - \zonesix ))\nad \zonefour &0\cr  - \zonefour &0& - (\zonefive  + \zonesix )&0\cr (\zonesix \cdot (2\cdot
\zonefive  + \zonesix ))\nad \zonefour &\zonefive  + \zonesix &0&0\cr 0&0&0&0}\koniec
,
\pmatrix{0&0&0&\zonefive \cr 0&0& - \zonefive &0\cr 0& - \zonefive &0&0\cr \zonefive &0&0&0}\koniec
$$
Here if $\zonesix  \neq 0$ we can make it vanish by the following symmetry.
\footnote{One may ask what happens when both $\zonefour$ and $\zonesix $ tend to $0$.
Then the proposed symmetry becomes singular and thus unapplicable.
We should notice that after performing that limit we would obtain
$$r_B=
\pmatrix{0&0&z&0\cr  0&0& - \zonefive   &0\cr 
 z &\zonefive  &0&0\cr 0&0&0&0}\koniec
$$
which differs from the case when $\zonesix =0$ with just the role of $X_+$
and $X_-$ interchanged and can be obtained as well by the symmetry
$\pmatrix{a & b\cr c & d}=\pmatrix{0 & 1\cr 1 & 0}$
}
$$\pmatrix{a& b \cr c & d} =
 \pmatrix{ 1 & - \frac{2\,\zonefive +\zonesix }{\zonefour } \cr 0  &1} $$
So we can assume that $\zonesix =0$:
$$
r_{\k{3}}=\pmatrix{0&\zonefour & 0&0\cr 
 - \zonefour &0& - \zonefive  &0\cr
  0 &\zonefive  &0&0
 \cr 0&0&0&0}
,
\pmatrix{0&0&0&\zonefive \cr 0&0& - \zonefive &0\cr 0& - \zonefive &0&0\cr \zonefive &0&0&0}\koniec
$$
Now if $\zonefive  \neq 0$ then we can make $\zonefour =0$ by using: 
$$\pmatrix{a& b \cr c & d} =
\pmatrix{0 &{( - 2\cdot d\cdot \zonefive )\over \zonefour } \cr 
{\zonefour \over (2\cdot d\cdot \zonefive )} & d}\,.$$
And we obtain a standard $r$-matrix
\bel{b1}
r_{b1}= x (X_+ \wedge X_- - V_+ \wedge W_- + V_- \wedge W_ +)\,.
\ee
where $x=-\zonefive $. Otherwise  $\zonefive =0$ and we we use 
$$\pmatrix{a& b \cr c & d} =
\pmatrix{1/\sqrt{\zonefour }&0 \cr 
0 & \sqrt{\zonefour }}\,.
$$
to obtain non standard $r$-matrix 
\bel{b2}
r_{b2}=  H \wedge X_+\,.
\ee

\subsection{Case \k{4}}
\pocz {\k{4}}
$$
\ba{ll}
r_{\k{4}}=&\pmatrix{
0&(\zfour\cdot \zzonesix )\over \zzten&(\zzten\cdot(\zzonefive ^2 - \zzonesix ^2))\over(\zfour\cdot \zzonesix )&\zzten
\cr  
( - \zfour\cdot \zzonesix )\over \zzten&0& - \zzonesix &\zfour\over2
\cr  
(\zzten\cdot( - \zzonefive ^2 + \zzonesix ^2))\over(\zfour\cdot \zzonesix )&\zzonesix &0&(\zzten^2\cdot( - \zzonefive ^2 + \zzonesix ^2))\over
(2\cdot \zfour\cdot \zzonesix ^2)
\cr  
 - \zzten&( - \zfour)\over2&(\zzten^2\cdot(\zzonefive ^2 - \zzonesix ^2))\over(2\cdot \zfour\cdot \zzonesix ^2)&0
}\\&
\pmatrix{
0&0&0&\zzonefive 
\cr  
0&0& - \zzonefive &0
\cr  
0& - \zzonefive &0&0
\cr  
\zzonefive &0&0&0
}
\ea
$$

In this case it is possible
\footnote{We assume that denominators of the entries should be different from 0.
Let us however look at the limits when both numerator and denominator tends to 0.
Please note that if $\zzonefive =0$ then $\zzonesix $ cancels from all the denominators.
We then obtain the purely bosonic r-matrix
$$
r_{\k{4}}=\pmatrix{
0&(\zfour\cdot \zzonesix )\over \zzten&(-\zzten\cdot \zzonesix )\over \zfour&\zzten
\cr  
( - \zfour\cdot \zzonesix )\over \zzten&0& - \zzonesix &\zfour\over2
\cr  
(\zzten\cdot \zzonesix )\over \zfour &\zzonesix &0&\zzten^2\over
(2\cdot \zfour)
\cr  
 - \zzten&( - \zfour)\over2&\zzten^2\over(2\cdot \zfour)&0
}
$$
Now, if $\zzonesix =0$ we obtain
$$
r_{\k{4}}=\pmatrix{
0&0&0&\zzten
\cr  
0&0& 0&\zfour\over2
\cr  
0 & 0&0&\zzten^2\over
(2\cdot \zfour)
\cr  
 - \zzten&( - \zfour)\over2&\zzten^2\over(2\cdot \zfour)&0
}
$$
which is always equivalent to 
$$
r_{\k{4}}=B\wedge X_+\,.
$$
If however $\zzonesix \neq 0$ then we may consider the limit $\zfour=\lambda \zzten \to 0$
and obtain
$$
r_{\k{4}}=\pmatrix{
0&\lambda\cdot \zzonesix &-\zzonesix \over \lambda &0
\cr  
 - \lambda\cdot \zzonesix &0& - \zzonesix &0 
\cr  
\zzonesix \over \lambda &\zzonesix &0& 0
\cr  
0&0&0&0
}
$$
which is equivalent to 
$$
r_{\k{4}}=H\wedge X_+\,.
$$
When $\zfour\neq0\neq \zzten$ we obtain ( symmetry 
$\pmatrix{a &b \cr c & d} = \pmatrix{1 &-d\cdot \zzten \over \zfour \cr
\zfour \cdot (1-d) \over d \cdot \zzten & d}$ with $d=\sqrt{-\zfour\over 2}$) 
$$
r_{c0} = x H\wedge X_+ + B \wedge X_+.
$$
The situation is more complicated when $\zzonefive \neq 0$. Then we have just 
two possibilities of singular limits:\
(a) $\zzonesix \to 0 $, $\zzten\to 0$, ${\zzten\over \zzonesix  }=$const.$=\zzten'$  
(b) $\zfour\to 0 $, $\zzten\to 0$, ${\zfour\over \zzten} =$const.$=\zfour'$.  
in case (a) we obtain:
$$
r_{\k{4}}=
\pmatrix{
0&\zfour\over \zzten'&\zzten'\cdot \zzonefive ^2 \over \zfour &0
\cr  
 - \zfour\over \zzten'&0& 0&\zfour\over2
\cr  
-\zzten'\cdot \zzonefive ^2\over \zfour &0&0&-\zzten'^2\cdot \zzonefive ^2 \over
(2\cdot \zfour)
\cr  
 0& - \zfour\over2& \zzten'^2\cdot \zzonefive ^2 \over 2\cdot \zfour&0
}\\
\pmatrix{
0&0&0&\zzonefive 
\cr  
0&0& - \zzonefive &0
\cr  
0& - \zzonefive &0&0
\cr  
\zzonefive &0&0&0
}
$$
which is equivalent to \r{c1}.

In case (b) we obtain:
$$
r_{\k{4}}=\pmatrix{
0&\zfour'\cdot \zzonesix &(\zzonefive ^2 - \zzonesix ^2)\over(\zfour'\cdot \zzonesix )&0
\cr  
 - \zfour'\cdot \zzonesix  &0& - \zzonesix & 0
\cr  
( - \zzonefive ^2 + \zzonesix ^2)\over(\zfour'\cdot \zzonesix )&\zzonesix &0&
0
\cr  
 0&0&0&0
},
\pmatrix{
0&0&0&\zzonefive 
\cr  
0&0& - \zzonefive &0
\cr  
0& - \zzonefive &0&0
\cr  
\zzonefive &0&0&0
}
$$
which is equivalent to \r{c1} with $y=0$.
}
 to perform the following symmetry:
$$\pmatrix{a& b \cr c & d} =
 \pmatrix{ 1 &   - { (\zzonesix +\zzonefive )\cdot \zzten\over \zfour\cdot \zzonesix  }\cr
  0 &  1} $$
and also rescale the parameter $\zzten$  as  $\zzten' = \zzten / \zzonesix  $ then the efect is
the same as setting $\zzonesix  =-\zzonefive $; we come up with the bosonic part:
$$
r_B=\pmatrix{
0&( - \zfour\cdot \zzonefive )\over \zzten'&0& \zzten'\cr
(\zfour\cdot \zzonefive )\over \zzten'&0&\zzonefive &\zfour\over2\cr
0& - \zzonefive &0&0\cr  - \zzten'&( - \zfour)\over2&0&0}
$$
and the $r_F$ unchanged. Then the symmetry
$$\pmatrix{a& b \cr c & d} =
\pmatrix{ 1 & 0 \cr  - {\zfour \over 2\cdot \zzten'} &1} 
$$
gives us effect of setting $\zfour=0$. We obtain a 2-parameter family
\bel{c1}
r_{c1}= y H \wedge B + x( X_+ \wedge X_- - V_+ \wedge W_- + V_-\wedge W_+)\,,
\ee
where $x=\zzonefive $ and $y=\zzten'$. 

\subsection{Cases \k{6},\k{11},\k{9}}
\pocz {\k{6}}
$$
r_{\k{6}}=\pmatrix{0&0& - 2\cdot \ztwo&0\cr 0&0&\zten&0\cr 2\cdot \ztwo& - \zten&0&0\cr 0&0&0&0}\koniec ,
\pmatrix{( - \zten\cdot \zonethree )\nad \ztwo&\zonethree &0&0\cr \zonethree &( - \zonethree \cdot \ztwo)\nad \zten&\zten& - \ztwo\cr 0&\zten&0&0\cr 0& -
\ztwo&0&( - \zten\cdot \ztwo)\nad \zonethree }\koniec
$$
Since $\ztwo$, $\zten$ and $\zonethree $ cannot vanish
\footnote{they cannot vanish seperately but any two of them can tend
to 0 at the same time. If $\ztwo$ and $\zten$ go to 0 we obtain a case
which is equivalent to $V_+\wedge V_-$. The remaining cases are equaivalent 
to  the one described below.
} we can perform the 
symmetry transformation 
$$\pmatrix{a & b \cr c &d } =
\pmatrix{ - {  \ztwo \over \sqrt{\ztwo\over \zonethree }\cdot \zten } & \sqrt{\ztwo \over \zonethree } \cr
  - {\sqrt{\zonethree \over \ztwo}} & 0 }
$$
we obtain one parameter family
$$
\t r_{\k{6}} =\pmatrix{0&0&0&0\cr 0&0& - \zten&0\cr 0&\zten&0&0\cr 0&0&0&0},
\pmatrix{0&0&0& - \zten\cr 0& - \zten&0&0\cr 0&0& - \zten&0\cr  - \zten&0&0&0}\,,
$$
which can be written as 
\bel{d1}
r_{d1} = x ( X_+ \wedge X_- + V_+ \wedge W_- + \frac12 V_- \wedge V_- +
\frac12 W _+ \wedge W_+ )\,,
\ee
where $x= - \zten$. The following two cases are equivalent to the above.   
\pocz {\k{11}}
$$
r_{\k{11}}=\pmatrix{0&0& - 2\cdot \ztwo&0\cr 0&0&\zonefive &0\cr 2\cdot \ztwo& - \zonefive &0&0\cr 0&0&0&0}\koniec ,
\pmatrix{0&0&0&\zonefive \cr 0&( - \zonefive \cdot \ztwo)\nad \zeight&0& - \ztwo\cr 0&0&( - \zonefive \cdot \zeight)\nad \ztwo&\zeight\cr \zonefive & - \ztwo&\zeight
&( - \ztwo\cdot \zeight)\nad \zonefive }\koniec
$$
This is equivalent to the Case \k{6} after application of $S$:
\pocz {\k{9}}
$$
r_{\k{9}}=\pmatrix{0& - 2\cdot \zone&0&0\cr 2\cdot \zone&0&\zten&0\cr 0& - \zten&0&0\cr 0&0&0&0}\koniec ,
\pmatrix{(\zone\cdot \zten)\nad \zeight&0&\zone&0\cr 0&0&\zten&0\cr \zone&\zten&(\zone\cdot \zeight)\nad \zten&\zeight\cr 0&0&\zeight&(\zten\cdot \zeight)\nad \zone
}\koniec
$$
\subsection{Cases \k{8},\k{22}}
\pocz {\k{8}}
$$
r_{\k{8}}=\pmatrix{0& - \zone& - \ztwo& - \zsix\cr \zone&0&\zsix\nad 2& - \zone\cr \ztwo&( - \zsix)\nad 2&0&\ztwo\cr \zsix&\zone& - \ztwo&0
}\koniec
,
\pmatrix{0&0&\zone&\zsix\nad 2\cr 0&0&\zsix\nad 2& - \ztwo\cr \zone&\zsix\nad 2&\zseven&\zeight\cr \zsix\nad 2& - \ztwo&\zeight&\zfive}\koniec
$$
\pocz {\k{22}}
$$
r_{\k{22}}=\pmatrix{0& - \zone& - \ztwo&\zonefive \cr \zone&0&\zonefive \nad 2&\zone\cr \ztwo&( - \zonefive )\nad 2&0& - \ztwo\cr  - \zonefive & - \zone&\ztwo
&0}\koniec
,
\pmatrix{\ztwelve&\zonethree &\zone&\zonefive \nad 2\cr \zonethree &\znine&\zonefive \nad 2& - \ztwo\cr \zone&\zonefive \nad 2&0&0\cr \zonefive \nad 2& - \ztwo&0&0}\koniec
$$
The case \k{8} is equivalent to case \k{22} upon the symmetry $S$ so we just focus on
the case \k{22}. 
We act in line with our general strategy. Since $r_{WW}=0$ we look at $r_{VW}$.
There are three  possibilities in this case\\
(a) $\rank (r_{VW}) = 2$,\\
(b) $\rank (r_{VW}) = 1$,\\
(c) $r_{VW}=0,\qquad (\zone=\ztwo=\zonefive =0)$.\\
In case (a) we can obtain a new matrix with the same structure but 
$\zone'=\ztwo'=0$ and $\zonefive ' = \sqrt{ \zonefive  ^2 - 4\cdot \zone \cdot \ztwo } \neq 0$. 
The symmetries preserving this form of $r_{VW}$ are generated by $(+-)$ and $(WV)$ scalings
and $(+-)$ swapping. By using this operations we can bring $r_{VV}$ to one 
of the following forms:
$$
\pmatrix{ \ztwelve' & \zonethree  ' \cr \zonethree ' &\znine' } =  
\pmatrix{ y & 1 \cr 1 & 1 } or   
\pmatrix{ 1 & 1 \cr 1 & 0 } or   
\pmatrix{ 1 & 0 \cr 0 & 1 } or   
\pmatrix{ \alpha & 0 \cr 0 & 0 } 
$$ 
This gives the following $r$-matrices
\bel{e0}
r_{e0} = x(2 H \wedge B + X_+ \wedge X_- + V_+ \wedge W_- + V_- \wedge W_+)\,,
\ee
\bel{e1}
r_{e1} = r_{e0}
         + y (V_+ \wedge V_+) + (V_+ \wedge V_-) +\frac12( V_- \wedge V_-) \,,
\ee
\bel{e2}
r_{e2} = r_{e0} + \frac12(V_+ \wedge V_+) + (V_+ \wedge V_-)  \,,  
\ee
\bel{e3}
r_{e3} = r_{e0} + \frac12 (V_+ \wedge V_+) + \frac12 (V_- \wedge V_-) \,,
\ee
\bel{e4}
r_{e4} = r_{e0}   +  \frac12 (V_+ \wedge V_+)\,,
\ee

In case (b) we can obtain the same structure with $\zone'=1$ and $\zonefive '=\ztwo'=0$.
i.e. $$r_{VW} = \pmatrix{ 1 &0 \cr 0 & 0}.$$
Symmetries preserving this form of $r_{VW}$ consists of $(WV)$ scaling 
combined with
$$
\pmatrix{a & b \cr c & d}= 
\pmatrix{1 & 0 \cr c & 1} 
$$ 
If $c9'\neq 0$ we can make $\zonethree '=0$, and then scale $\znine'$ 
to become 1. If $\znine'=0$ but $\zonethree '\neq 0$ we use c to make $\ztwelve'$ vanish.
If $\znine'=\zonethree '=0$ we just scale $\ztwelve'$ to 1 if is is not zero.
Thus the following possibilities for $r_{VV}$ emerge
$$
\pmatrix{ \ztwelve' & \zonethree  ' \cr \zonethree ' &\znine' } =  
\pmatrix{ y & 0 \cr 0 & 1 } or   
\pmatrix{ 0 & 1 \cr 1 & 0 } or   
\pmatrix{ \alpha & 0 \cr 0 & 0 }\,. 
$$ So we have
\bel{e5}
r_{e5}= - (H + B) \wedge X_+ + V_+ \wedge W_+ 
\ee
\bel{e6}
r_{e6}= r_{e5} + y \frac12(V_+ \wedge V_+ ) + \frac12( V_- \wedge V_-)  \,,
\ee
\bel{e7}
r_{e7}= r_{e5} + (V_+ \wedge V- )\,,    
\ee
\bel{e8}
r_{e8} = r_{e5} + \frac12(V_+ \wedge V+ )\,.
\ee
In case (c) we have all the symmetry and we can obtain 
depending on the rank of $r_{VV}$  
$$
\pmatrix{ \ztwelve' & \zonethree  ' \cr \zonethree ' &\znine } =  
\pmatrix{ 0 & 1 \cr 1 & 0 } or   
\pmatrix{ \alpha & 0 \cr 0 & 0 } 
$$ 
\bel{e9}
r_{e9} =  (V_+ \wedge V_- )\,.
\ee
\bel{e10}
r_{e10} =  \frac12(V_+ \wedge V_+ )\,.
\ee

\subsection{Case \k{10}}
\pocz {\k{10}}
$$
r_{\k{10}}=\pmatrix{0& - \zone& - \ztwo&0\cr \zone&0&0&\zfour\nad 2\cr \ztwo&0&0&( - \ztwo\cdot \zfour)\nad (2\cdot \zone)\cr 0&( - \zfour)\nad 2&(\ztwo\cdot
\zfour)\nad (2\cdot \zone)&0}\koniec
,
\pmatrix{0&0&\zone&0\cr 0&0&0& - \ztwo\cr \zone&0&0&0\cr 0& - \ztwo&0&0}\koniec
$$
When $\ztwo=0$ then $\zone$ may also vanish. In such a case 
after a simple ($+-$) rescaling we obtain just 
\bel{f0}
r_{f0}= B\wedge X_+
\ee
If $\ztwo=0$ but $\zone\neq0$ we just take 
\footnote{if $\zone$ and $\ztwo$ tend to 0 in such a way that
$ \frac{\ztwo}{\zone} $ 
remains finite we are left with a linear combination of $B\wedge X_+$
and $B\wedge X_-$ which is equivalent either to $B\wedge X_+$ or $y H \wedge B$.
}
$$\pmatrix{a& b \cr c & d} =
\pmatrix{1/\sqrt{\zone} &0 \cr 0 &\sqrt{\zone}}
$$ 
and obtain
$$
\t r_{\k{10}}=\pmatrix{0& - 1& 0&0\cr 1&0&0&\zfour\nad 2\cr 0&0&0& 0\cr 
 0&( - \zfour)\nad 2&0&0}\koniec
,
\pmatrix{0&0&1&0\cr 0&0&0& 0\cr 1&0&0&0\cr 0& 0&0&0}\koniec
$$
which can be written as
\bel{f1}
r_{f1}= (-H\wedge X_+ + V_+ \wedge W_+) + x (B \wedge X_+)\,,
\ee
where $x=-\frac {\zfour}{2}$. If, however, $\ztwo\neq 0$ than we take 
$$
\pmatrix{a& b \cr c & d} =
\pmatrix{1/\sqrt{\zone} &-\frac12\sqrt{\ztwo}\cr  1/\sqrt{\ztwo} &\frac12\sqrt{\zone} }
$$
to obtain
$$\t{\t r}_{10}=
\pmatrix{0&0&0& y\cr 
0& 0&x&0\cr 
0&-x& 0&0\cr  
-y&0&0&0}{{}}
,
\pmatrix{0&0&0& x\cr 
0& 0&x&0\cr 
0&x& 0&0\cr  
x&0&0&0}{{}}
$$
or 
\bel{f2}
r_{f2} = x( X_= \wedge X_- + V_+ \wedge W_- + V_- \wedge W+) 
+ y(H\wedge B) \,,
\ee
where $x=-\sqrt{\zone\cdot \ztwo} $ and $y=-\zfour \sqrt{\ztwo\over \zone}$.
\subsection{Cases \k{12},\k{18}}
\pocz {\k{12}}
$$
r_{\k{12}}=\pmatrix{0&0&( - \zthree)\nad 2&\zonefive \cr 0&0&\zonefive \nad 2&0\cr \zthree\nad 2&( - \zonefive )\nad 2&0&( - \zthree)\nad 2\cr  - \zonefive &0&\zthree
\nad 2&0}\koniec ,
\pmatrix{0&0&0&\zonefive \nad 2\cr 0&0&( - \zonefive )\nad 2&0\cr 0&( - \zonefive )\nad 2&( - 2\cdot \zonefive \cdot \zeight)\nad \zthree&\zeight\cr \zonefive \nad 2&0
&\zeight&( - \zthree\cdot \zeight)\nad (2\cdot \zonefive )}\koniec
$$
Here if $\zeight=0$ and $\zonefive \neq0$ then we take 
$$
\pmatrix{a& b \cr c & d} =
\pmatrix{1 &\frac{\zthree}{2\cdot \zonefive }\cr 0 &1}
$$ to obtain
$$
\t r_{\k{12}}=\pmatrix{0&0&0&\zonefive \cr 0&0&\zonefive \nad 2&0\cr 0&( - \zonefive )\nad 2&0&0\cr
  - \zonefive &0& 0&0}\koniec ,
\pmatrix{0&0&0&\zonefive \nad 2\cr 0&0&( - \zonefive )\nad 2&0\cr
 0&( - \zonefive )\nad 2&0 & 0\cr 
 \zonefive \nad 2&0 & 0&0}\koniec
$$
This is just 
\bel{g1}
r_{g1} = x( 2 H\wedge B + X_+ \wedge X_- + V_+ - V_- \wedge W_+).
\ee 
with $x=\frac{\zonefive }{2}$.
If $\zeight=0$ then we can also put $\zonefive =0$ and obtain after rescaling 
\bel{g1.5}
r_{g1.5} =(H-B) \wedge X_+\,.
\ee
If $\zeight\neq 0$ we take 
\footnote{if $\zthree$ and $\zonefive $ tend to 0 we are left with a purely fermionic 
$r$-matrix which is easily shown to be equaivalent with $V_+\wedge V_+$.}
$$\pmatrix{a& b \cr c & d} =
\pmatrix{\frac{\ztwo}{2\cdot b \cdot \zonefive } &\sqrt{-\zthree \cdot \zeight \over 2 
\cdot \zonefive } \cr 
0 &  \frac{2\cdot b \cdot \zonefive }{\ztwo}}
$$ and 
after performing additional $S$ tranformation obtain
$$
\t{\t r}_{12}=\pmatrix{0&0&0&\zonefive \cr 0&0&\zonefive \nad 2&0\cr 0&( - \zonefive )\nad 2&0&0\cr
  - \zonefive &0& 0&0}\koniec ,
\pmatrix{1&0&0&\zonefive \nad 2\cr 0&0&( - \zonefive )\nad 2&0\cr
 0&( - \zonefive )\nad 2&0 & 0\cr 
 \zonefive \nad 2&0 & 0&0}, 
$$
which differs from \r{g1} only by $\frac12 V_+ \wedge V_+$
\bel{g2}
r_{g1} = r_{g1} + \frac12 V_+ \wedge V_+\,.
\ee
\pocz {\k{18}} 
$$
r_{\k{18}}=\pmatrix{0&0&\zthree\nad 2&\zonefive \cr 0&0&( - \zonefive )\nad 2&0\cr ( - \zthree)\nad 2&\zonefive \nad 2&0&( - \zthree)\nad 2\cr  - \zonefive &0&\zthree
\nad 2&0}\koniec ,
\pmatrix{( - 2\cdot \zonethree \cdot \zonefive )\nad \zthree&\zonethree &0&\zonefive \nad 2\cr
\zonethree &( - \zonethree \cdot \zthree)\nad (2\cdot \zonefive  ) & (  - \zonefive )\nad 2&0\cr
0&( - \zonefive )\nad 2&0&0\cr
\zonefive \nad 2&0&0&0}\koniec
$$
This case  is equivalent to the Case \k{12}  upon the symmetry $S$.
\subsection{Cases \k{20},\k{2}}
The matrices $r_{\k{2}}$ and $r_{\k{20}}$ differ only by the names of parameters \\
$$
\ba{ll}
r_{\k{20}}=&\pmatrix{0&0&(\zthree\cdot ( - \zten + \zonefive ))\nad (2\cdot (\zten + \zonefive ))&\zten + \zonefive \cr
 0&0&(\zten - \zonefive )\nad 2&0\cr
 ( \zthree\cdot (\zten - \zonefive ))\nad (2\cdot (\zten + \zonefive ) ) & (  - \zten + \zonefive )\nad 2&0&( - \zthree)\nad 2\cr  - (\zten + \zonefive )&0&
\zthree\nad 2&0}
\\[12mm]
&\pmatrix{0&0&0&( - \zten + \zonefive )\nad 2\cr 0&0&(\zten - \zonefive )\nad 2&0\cr 0&(\zten - \zonefive )\nad 2&0&0\cr ( -
\zten + \zonefive )\nad 2&0&0&0}
\ea
$$
$$
r_{\k{2}}=\pmatrix{0&0&(\zonesix \cdot \zthree)\nad (2\cdot (2\cdot \zonefive  + \zonesix ))&2\cdot \zonefive  + \zonesix \cr 0&0&( - \zonesix )\nad 2&0\cr ( - \zonesix \cdot \zthree)
\nad (2\cdot (2\cdot \zonefive  + \zonesix ))&\zonesix \nad 2&0&( - \zthree)\nad 2\cr  - 2\cdot \zonefive  - \zonesix &0&\zthree\nad 2&0}\koniec ,
\pmatrix{0&0&0&( - \zonesix )\nad 2\cr 0&0&\zonesix \nad 2&0\cr 0&\zonesix \nad 2&0&0\cr ( - \zonesix )\nad 2&0&0&0}\koniec
$$
so just look into the $r_{\k{2}}$. The symmetry given by 
\footnote{If $\zonefive $ and $\zonesix $ both tend to 0 we obtain any linear combination
of $H\wedge X_-$ and $B\wedge X_-$.
}
$$\pmatrix{a& b \cr c & d} =
\pmatrix{1 &{ \zthree \over 2 (2 \cdot \zonefive  + \zonesix )} \cr 0  & 1} 
$$ 
yields 
$$
\tilde r_{\k{2}}=\pmatrix{0&0&0&2\cdot \zonefive  + \zonesix \cr 
0&0&( - \zonesix )\nad 2&0\cr 
0&\zonesix \nad 2&0& 0\cr
  - 2\cdot \zonefive  - \zonesix &0&0&0}\koniec ,
\pmatrix{0&0&0&( - \zonesix )\nad 2\cr 0&0&\zonesix \nad 2&0\cr 0&\zonesix \nad 2&0&0\cr ( - \zonesix )\nad 2&0&0&0}\,,
$$
which can be written as
\bel{h1}
r_{h1} = x( X_+ \wedge X_- + V_+ \wedge W_- - V_- \wedge W_+)
+ y( H \wedge B)\,,
\ee
with $x=\frac{-\zonesix }{2}$ and $y=2 \cdot \zonefive  + \zonesix $.
\subsection{Cases \k{14},\k{15}}
\pocz {\k{14}}
$$
r_{\k{14}}=\pmatrix{0&(\sqrt{\zfour}\cdot \ztwoone)\nad \sqrt{\zthree} & 
( \sqrt{\zfour}\cdot \ztwoone\cdot \zthree)\nad (\sqrt{\zthree}\cdot \zfour)&0\cr 
( - \sqrt{\zfour}\cdot
\ztwoone)\nad \sqrt{\zthree}&0&0&\zfour\nad 2\cr 
( - \sqrt{\zfour}\cdot \ztwoone\cdot \zthree)\nad (\sqrt{\zthree}\cdot \zfour)&0&0&( - \zthree)\nad 2\cr 0&
(- \zfour)\nad 2&\zthree\nad 2&0}\koniec
,
\pmatrix{0&0&0&\ztwoone\cr 0&0& - \ztwoone&0\cr 0& - \ztwoone&0&0\cr \ztwoone&0&0&0}\koniec
$$
From the above we see that both $\zthree$ and $\zfour$ are non 
zero \footnote{However, we can look at the limit when both $\zthree$  and $\zfour$ tend to $0$,
leaving the ratio $\zthree\over \zfour$ finite. 
The same symmetry is applicable this case, leading to the
$r$-matrix $r_{i1}$ with $y=0$.}.
Upon the symmetry 
$$\pmatrix{a& b \cr c & d} =
\pmatrix{ \frac12 & \sqrt{\zthree \over \zfour}\cr 
 -\frac12 \sqrt{\zfour \over \zthree}& 1}
$$ 
 we obtain
$$
\tilde r_{\k{14}}=
\pmatrix{0&0&0&\sqrt{\zthree\cdot \zfour}
 \cr 0&0& - \ztwoone&0\cr
  0&  \ztwoone&0&0\cr 
  -\sqrt{\zthree\cdot \zfour}
  &0&0&0}\koniec
,
\pmatrix{0&0&0&\ztwoone\cr 0&0& - \ztwoone&0\cr 0& - \ztwoone&0&0\cr \ztwoone&0&0&0}.
$$ 
This can be written as 
\bel{i1} 
r_{i1} = x (X_+ \wedge X_- - V_+ \wedge W_- + V_- \wedge W_+)
         + y (H\wedge B)
\ee
with $x = -\ztwoone$ and $y= \sqrt{\zthree \cdot \zfour}$.

If $\ztwoone=0$ we have 
$$
r_{\k{14}}=\pmatrix{0&0& 0&0\cr 
0&0&0&\zfour\nad 2\cr 
0&0&0&( - \zthree)\nad 2\cr 0&
(- \zfour)\nad 2&\zthree\nad 2&0}\koniec
,
\pmatrix{0&0&0&0\cr 0&0& 0&0\cr 0& 0&0&0\cr 0&0&0&0}\koniec
$$
then it is possible that just one of $\zthree$ and $\zfour$ vanishes.
After putting $\zthree=0$ and rescaling we obtain
\bel{i2}
r_{i2}= B \wedge X_+\,.
\ee
Case when $\ztwoone=0$ and $\zfour\neq0\neq \zthree$ is just $r_{i1}$ with $x=0$.

\pocz {\k{15}} is equivalent to $r_{\k{14}}$ the only difference being the sign of
$\sqrt{\zfour}$ and $\ztwoone$ renamed $\ztwotwo$.
$$
r_{\k{15}}=\pmatrix{0&( - \sqrt{\zfour}\cdot \ztwotwo)\nad \sqrt{\zthree} & 
(  - \sqrt{\zfour}\cdot \ztwotwo\cdot \zthree)\nad (\sqrt{\zthree}\cdot \zfour)&0\cr 
(\sqrt{\zfour}\cdot \ztwotwo)\nad \sqrt{\zthree}&0&0&\zfour\nad 2\cr 
(\sqrt{\zfour}\cdot \ztwotwo\cdot \zthree)\nad (\sqrt{\zthree}\cdot \zfour)&0&0&( - \zthree)\nad 2\cr 0&(
 - \zfour)\nad 2&\zthree\nad 2&0}\koniec
,
\pmatrix{0&0&0&\ztwotwo\cr 0&0& - \ztwotwo&0\cr 0& - \ztwotwo&0&0\cr \ztwotwo&0&0&0}\koniec
$$
\subsection{Case \k{17}}
\pocz {\k{17}}
$$
r_{\k{17}}=\pmatrix{0& - 2\cdot \zone&(2\cdot \zonethree \cdot \zeight)\nad \zone&0\cr 2\cdot \zone&0&(\zonethree \cdot \zeight + \zonefive ^2)\nad \zonefive &0\cr ( - 2\cdot \zonethree \cdot \zeight)\nad
\zone&( - (\zonethree \cdot \zeight + \zonefive ^2))\nad \zonefive &0&0\cr 0&0&0&0}\koniec
,
\pmatrix{(\zone\cdot \zonethree )\nad \zonefive &\zonethree &\zone&\zonefive \cr \zonethree &(\zonethree \cdot \zonefive )\nad \zone&(\zonethree \cdot \zeight)\nad \zonefive &(\zonethree \cdot \zeight)\nad \zone\cr \zone&(
\zonethree \cdot \zeight)\nad \zonefive &(\zone\cdot \zonefive )\nad \zonethree &\zeight\cr \zonefive &(\zonethree \cdot \zeight)\nad \zone&\zeight&(\zonethree \cdot \zeight^2)\nad (\zone\cdot \zonefive )}\koniec
$$
If $\zeight=0$ we have\footnote{The singular limits were investigated but gave rise
to no new $r$-matrices. We therefore assume that 
all denominators are different from 0.}
$$
r_{\k{17}}=\pmatrix{0& - 2\cdot \zone&0&0\cr 2\cdot \zone&0&\zonefive &0\cr 
0&- \zonefive &0&0\cr 0&0&0&0}\koniec
,
\pmatrix{(\zone\cdot \zonethree )\nad \zonefive &\zonethree &\zone&\zonefive \cr \zonethree &(\zonethree \cdot \zonefive )\nad \zone&0& 0\cr \zone&
0&(\zone\cdot \zonefive )\nad \zonethree &0\cr \zonefive &0&0&0}\koniec
$$
then after symmetry with 
$$
\pmatrix{a& b \cr c & d} =
\pmatrix{\sqrt{\zonethree \over \zone} & 0 \cr
-\sqrt{\zonethree \cdot \zone}/\zonefive  & \sqrt{\zone \over \zonethree }}
$$ we obtain
$$
r_{\k{17}}=\pmatrix{0&0&0&0\cr 0&0&\zonefive &0\cr 
0&- \zonefive &0&0\cr 0&0&0&0}\koniec
,
\pmatrix{ 0&0&0&\zonefive \cr 0&\zonefive &0& 0\cr 
0& 0&\zonefive &0\cr 
\zonefive &0&0&0}\koniec
$$
which is just 
\bel{j1}
r_{j1} = x( X_+ \wedge X_- + V_+ \wedge W_- + \frac12 V_- \wedge V_-
+ \frac12 W_+ \wedge  W_+ )\,,
\ee
with $x=\zonefive $.
When $\zeight \neq 0$ we still have two possibilities.
If $\zonefive ^2-\zonethree \cdot \zeight\neq 0$ we can obtain a similar matrix (with $\zonefive $ replaced with 
$ \zonefive ^2 -\zonethree  \cdot \zeight \over \zonefive $). To do it we take
$$\pmatrix{a& b \cr c & d} =
\pmatrix{-c\cdot \zonefive  / \zone & \zonethree \cdot \zeight \over c \cdot (\zonefive ^2-\zonethree \cdot \zeight) \cr
{ \zone\cdot \zonethree  \over \zonefive ^2 -\zonethree  \cdot \zeight } & 
{ -\zone\cdot \zonefive  \over c \cdot (\zonefive ^2-\zonethree \cdot \zeight) }}\,.
$$
If $\zonefive ^2-\zonethree \cdot \zeight = 0$ we then let $\zonethree =\zonefive ^2/\zeight$ 
and the symmetry
$$\pmatrix{a& b \cr c & d} =
\pmatrix{{ -1\over \sqrt{\zone}} &{\zonefive \over \sqrt{\zone}}\cr
0 & -\sqrt{\zone}}
$$ 
yields
$$
r_{\k{17}}=\pmatrix{0&-2&0&0\cr 2&0&0&0\cr 
0&0&0&0\cr 0&0&0&0}\koniec
,
\pmatrix{\zonefive  \nad \zeight &0&1&0\cr 
0&0&0& 0\cr 
1& 0&\zeight\nad \zonefive &0\cr 
0&0&0&0}\koniec
$$
When we take $k=\zonefive /\zeight$ and it becomes   
$$
r_{\k{17}}=\pmatrix{0&-2&0&0\cr 2&0&0&0\cr 
0&0&0&0\cr 0&0&0&0}\koniec
,
\pmatrix{1&0&1&0\cr 
0&0&0& 0\cr 
1& 0&1&0\cr 
0&0&0&0}\,,
$$
which can be written as
\bel{j2}
r_{j2} = -2 H\wedge X_+ + \frac12( V_+ + W_+) \wedge ( V_+ + W_+)\,.
\ee 
\subsection{Case \k{21}}
\pocz {\k{21}}
$$
\ba{ll}\dsp
r_{\k{21}}=&\pmatrix{0& - \zone& - \ztwo&\zonefive  - \zsix\cr
\zone&0&(\zonefive  + \zsix)\nad 2&(\zone\cdot (\zonefive  - \zsix))\nad (\zonefive  + \zsix)\cr
 \ztwo& ( - (\zonefive  + \zsix))\nad 2&0&(\ztwo\cdot ( - \zonefive  + \zsix))\nad (\zonefive  + \zsix)\cr
  - \zonefive  + \zsix&(\zone\cdot ( - \zonefive  + \zsix ))\nad (\zonefive  + \zsix) & ( \ztwo\cdot (\zonefive  - \zsix))\nad (\zonefive  + \zsix)&0}\koniec
\\
&\pmatrix{0&0&\zone&(\zonefive  + \zsix)\nad 2\cr 0&0&(\zonefive  + \zsix)\nad 2& - \ztwo\cr \zone&(\zonefive  + \zsix)\nad 2&0&0\cr (\zonefive  +
 \zsix)\nad 2& - \ztwo&0&0}\koniec
 \ea
$$
Let us denote $ \zonefive =\zzonefive +\zzsix$ $\zsix=\zzonefive -\zzsix$. Then we have two possibilities; 
\ben
\item[(a)] $\zone\cdot \ztwo+\zzonefive ^2=0$,
\item[(b)] $\zone\cdot \ztwo+\zzonefive ^2\neq 0$.
\een
In the latter case both $\zone$ and $\ztwo$ can be made equal zero.
If $\ztwo=0$ we just use the matrix 
$$
\pmatrix{a &b \cr c & d} =
\pmatrix{0 &-1 \cr 1 & \zone\over \zonefive  +\zsix } 
$$
and when $\ztwo\neq 0$ we use
$$
\pmatrix{a &b \cr c & d} =
\pmatrix{ \frac{\ztwo}{2\sqrt{\zone\cdot \ztwo + \zzonefive ^2}} &  1 \cr
\frac{\zonefive }{2\sqrt{\zone\cdot \ztwo + \zzonefive ^2}} - \frac12 &  
\frac{\zzonefive +\sqrt{\zone\cdot \ztwo + \zzonefive ^2}}{\ztwo}} $$
in either case we obtain  
$$
\ba{ll}
r_{\k{21}}=&\pmatrix{0& 0& 0&\zonefive  - \zsix\cr
0&0&(\zonefive  + \zsix)\nad 2&0\cr
 0& ( - (\zonefive  + \zsix))\nad 2&0&0\cr
  - \zonefive  + \zsix&  0 & 0&0}\koniec
\\
&\pmatrix{0&0&0&(\zonefive  + \zsix)\nad 2\cr 0&0&(\zonefive  + \zsix)\nad 2& 0\cr
0&(\zonefive  + \zsix)\nad 2&0&0\cr (\zonefive  +  \zsix)\nad 2& 0&0&0}\koniec
\ea
$$
This is just
\bel{k1}
r_{k1} = x( X_+ \wedge X_- + V_+ \wedge W_- + V_- \wedge W_+ ) + y (H\wedge B)
\ee
with $x=\frac{\zonefive +\zsix}{2}$ and $y = \zonefive -\zsix$.  
In the case (a), however,  we use the symmetry 
$$
\pmatrix{a &b \cr c & d} =
\pmatrix{0 &\sqrt{-\ztwo} \cr \sqrt{-\ztwo}\over \ztwo &
\sqrt{-\ztwo}(\zonefive  +\zsix)\over 2 \cdot \ztwo } 
$$
and obtain 
$$
r_{\k{21}}=\pmatrix{0& 1& 0&0\cr
-1&0&0 &\zonefive  -\zsix  \nad \zonefive  + \zsix\cr
 0& 0&0&0\cr
 0& - \zonefive  + \zsix\nad \zonefive +\zsix  & 0&0},
\pmatrix{0&0&1& 0\cr 
0&0&0& 0\cr
1&0&0&0\cr 
0& 0&0&0}
$$
or
\bel{k2}
r_{k2} = -( H \wedge X_+) + ( V_+ \wedge W_+) + x(B \wedge X_+)\,,
\ee  
where $x= \frac{\zonefive -\zsix }{\zonefive +\zsix}$.

\newpage
\section{Summary and discussion}
Below we give the list nonequivalent $r$-matrices  for the $osp(2|2)$.
The lowercase latin letters $x$, $y$ denote arbitrary complex numbers, 
whereas $\alpha$, $\beta$ can only take value 0 or 1. 
\begin{eqnarray}
\label{b2}
r_{b2}&=&  H \wedge X_+\,,~~~~~~~~~~~~~~~~~~~~~~~~~~~~~~~~~~~~~~~~~~~~~~~~~~~~~~~~~~~~~~~~~~~~~~~~~~~~~~~~
\\
\label{c0}
r_{c0} &=& x H\wedge X_+ + B \wedge X_+\,,
\\
\label{a2}
r_{a2} &=& \alpha (H - B) \wedge X_+ + \beta (V_+ \wedge V_+)\,, 
\\
\label{h1}
\hspace*{-4cm}r_{b1}\subset r_{c1}=r_{i1} \sim r_{h1} &=& x( X_+ \wedge X_- + V_+ \wedge W_- - V_- \wedge W_+)
+ y( H \wedge B)\,,
\\
\label{f2}
r_{f2} =r_{k1}&=& x( X_+ \wedge X_- + V_+ \wedge W_- + V_- \wedge W_+) + y(H\wedge B) \,,
\\
\label{d1}
r_{d1} =r_{j1}&=& x ( X_+ \wedge X_- + V_+ \wedge W_- + \frac12 V_- \wedge V_- +
\frac12 W _+ \wedge W_+ )\,,
\\
\label{j2}
r_{j2} &=& -2 (H\wedge X_+)  + \frac12( V_+ + W_+) \wedge ( V_+ + W_+)\,,
\\ 
\label{g}
r_{g} &=& x( 2 H\wedge B + X_+ \wedge X_- + V_+\wedge W_- - V_- \wedge W_+)
+ \alpha \frac12 (V_+ \wedge V_+)\,,
\\
\label{a1}
r_{a1} &=& x(- 2 H \wedge B + X_+ \wedge X_- + V_+ \wedge W_- - V_- \wedge W_+)
+ \alpha \frac12  (V_+ \wedge V_+) \,,
\\
\label{e0}
r_{e0} &=& x(2 H \wedge B + X_+ \wedge X_- + V_+ \wedge W_- + V_- \wedge W_+)\,,
\\
\label{e1}
r_{e1} &=& r_{e0}
         + y (V_+ \wedge V_+) + (V_+ \wedge V_-) + \frac12 (V_- \wedge V_-) \,,
\\
\label{e2}
r_{e2} &=& r_{e0} + \frac12(V_+ \wedge V_+) + (V_+ \wedge V_-)  \,,  
\\
\label{e3}
r_{e3} &=& r_{e0} + \frac12 (V_+ \wedge V_+) + \frac12 (V_- \wedge V_-) \,,
\\
\label{e4}
r_{e4} &=& r_{e0}   +  \frac12 (V_+ \wedge V_+)\,,
\\
\label{f1}
r_{f1}=r_{k2} &=& ((xB-H )\wedge X_+) + (V_+ \wedge W_+) \,,
\\
\label{e5}
r_{e5} &\equiv&  -(B + H) \wedge X_+ + (V_+ \wedge W_+)\,, 
\\
\label{e6}
r_{e6} &=& r_{e5} + y \frac12(V_+ \wedge V_+ ) + \frac12( V_- \wedge V_-)  \,,
\\
\label{e7}
r_{e7} &=& r_{e5} + (V_+ \wedge V- )\,,    
\\
\label{e8}
r_{e8} &=& r_{e5} + \frac12(V_+ \wedge V+ )\,,
\\
\label{e9}
r_{e9} &=&  V_+ \wedge V_- \,,
\\
\label{e10}
r_{e10} &=&  \frac12(V_+ \wedge V_+ )\,.
\end{eqnarray}


All the ``generic'' solutions were found by the computer.
We carefully analysed the singular limits to obtain some more solutions.
On this basis we claim that every co-Lie structure on the $osp(2|2)$ algebra
must be generated by an $r$-matrix which is equivalent to some member of the 
list \r{b2}-\r{e10}. Due to the equivalence 
$osp(2|2) \sim sl(1|2)$ (see e.g\ \cite{kac}) this classification is also valid for the
$sl(1|2)$ super-algebra.

The $r$-matrices \r{b2}, \r{c0}, \r{a2}, \r{j2}, \r{f1}, \r{e5}, 
\r{e6}, \r{e7}, \r{e8}, \r{e9}, \r{e10} satisfy CYBE. 
The remaining ones  satisfy CYBE only if the parameter $x$ is equal to 0.

If the result of Etingof \cite{geom} can be generalized
to the case of Lie super-bialgebras,
then $r$-matrices satisfying CYBE can be easily quantized.

In view of the sequence of inclusions:
\bel{inc}
sl(2) \subset {gl(2) \atop osp(1|2)} \subset osp(2|2)
\ee
it is relevant to look at the classification of each subalgebra in this 
chain. 
It is obvious that any nonstandard $r$-matrix 
(i.e.\ satisfying CYBE) of a subalgebra 
is also $r$-matrix for the whole algebra. 
We make sure that nonstandard $r$-matrices 
of  $sl(2)$, $gl(2)$ and $osp(1|2)$ which are known in the literature 
are also present in our classification.
\begin{enumerate}
\item[(a)] for $sl(2)$ all the $r$-matrices satisfying CYBE are equivalent to $H\wedge 
X_+$ \r{b2}.

\item[(b)] The classification of $gl(2)$ Lie bialgebras was first obtained 
by Ballesteros et al.\ in \cite{balles},
where also the corresponding Hopf algebras were described.
From their nonstandard $r$-matrices it is possible to pick up just two 
nonequivalent $r_1=H\wedge X_+$ and $r_2=H \wedge B$ 
which coincide with \r{b2} and \r{c0}. 

\item[(c)]
The subalgebra $osp(1|2)$ is generated by $H, X_+, X_- ,V'_+=(V_++W_+)/2$ 
and $V'_-=(V_-+W_-)/2$. 
The classification of super-Lie bialgebras was obtained in \cite{js2}.
There were two nonstandard $r$-matrices:
$r_1=H\wedge X_+$ and $r_2=H\wedge X_+ - V_+ \wedge V_+$. 
They correspond to \r{b2} and \r{j2} from our list.
\end{enumerate}

The classification of $osp(2|2)$ Lie superbialgebras was not
known before, however, several examples have been investigated.

Deguchi et.al.\ \cite{deguchi} constructed  the deformation 
of the universal enveloping algebra $U_q(osp(2|2))$ and obtained its 
universal $R$-matrix..
After identification of the generators $J_\pm=\pm X_\pm$, $V_\pm=V_\pm/\sqrt2$,
$\overline V_\pm=W_\pm/\sqrt2$, $H=H$, $T =2B$ 
we notice that  the antisymmetric part of the first order term 
(in $\ln q$) of the $R$-matrix
takes the form  
$X_+ \wedge X_- + V_+\wedge W_- + V_-\wedge W_+$ 
which is a special case of our classical $r$-matrix \r{f2}.
We also check that it generates the antisymmetric part of the 
first order term of their coproduct. 

The universal $R$-matrix given by Aizawa \cite{aizawa} 
was obtained by a twisting element belonging to the $gl(2)$ subalgebra 
which had the following form (we use the following identification of 
generators used in \cite{aizawa}: $H=2H$, $Z=2B$, $X_\pm=\pm X_\pm$,
$v_\pm=V_\pm$, $\overline v_\pm = \pm W_\pm$ in order to 
give the original expression in our basis):
\bel{twist}
{\cal F}=\exp( \frac{g}{h} \sigma \tens B ) \exp (-H \tens \sigma)
\ee
where
$$ \sigma \equiv -\ln(1-2h X_+)
$$ 
The universal $R$-matrix takes the form:
$$
R = \exp(\frac{g}{h} B \tens \sigma ) 
    \exp(-\sigma \tens H)
    \exp(H\tens \sigma)
    \exp(-\frac{g}{h}\sigma \tens B)
$$
and in the classical limit $h\to 0$, $g\to 0$
 it gives rise to the  $r$-matrix \r{c0}.

Another two parameter deformation was investigated by Arnaudon et al. 
\cite{frappat}. After the identification of the generators 
($E_1^+=V_-$, $E_1^-=W_+$, $E_2^+=W_-$, $E_2^-=V_+$, $E_3^+=V_-$,
$E_3^-=X_+$, $H_1=H+B$, $H_2=H-B$)  
we were able to check that the 
super antisymmetric part of the first order term of their coproduct 
is generated by the classical $r$-matrix \r{f2}.

\section{Classification of $osp(1|2) \oplus u(1) $ super Lie bialgebras}
The $osp(1|2) \oplus u(1) $ Lie superalgebra has the same subalgebra
structure as $osp(2|2)$.
\bel{inc}
sl(2) \subset {gl(2) \atop osp(1|2)} \subset osp(1|2) \oplus u(1) 
\ee
However, it has only 6 generators so it is relatively 
easy to classify using the same technique.
 The $osp(1|2)$ algebra of 
is spanned by the generators $H$, $X_+$, $X_-$,
$Q_+\equiv \frac12(V_++W_+)$ and $Q_-\equiv \frac12(V_-+W_-)$, whose 
commutation relations follow from \r{r1}.
Supplementing them with a central generator $Z$ gives $osp(1|2) \oplus u(1)$.

All the Lie super-bialgebras $osp(1|2) \oplus u(1) $ are coboundary and 
their corresponding $r$-matrices are equivalent to one of the following
\beq
\label{o1}
r_1&=&H \wedge X_+\,,\\
\label{o1a}
r_2&=&Z \wedge X_+\,,\\
\label{o2}
r_3&=&H \wedge X_+ + Z \wedge X_+\,,\\
\label{o3}
r_4&=& H\wedge X_+ -  Q_+ \wedge Q_+ \,,\\ 
\label{o4}
r_5&=& H\wedge X_+ -  Q_+ \wedge Q_+  + Z \wedge X_+\,,\\ 
\label{o5}
r_6 &=& x(X_+ \wedge X_- + 2 Q_+ \wedge Q_- ) \,, \\
\label{o6}
r_7 &=& x(X_+ \wedge X_- + 2 Q_+ \wedge Q_- ) + H \wedge Z\,. 
\eeq
$r$-matrices \r{o1}-\r{o4} satisfy CYBE whereas \r{o5} and \r{o6} don't
if $x\neq 0$.

\subsection*{Acknowledgement}
The work on this paper was supported by the KBN grant No 2PO31313012.
The author would like to thank Prof.\ J.T.\ Sobczyk for stimulating discussions.

\end{document}